\documentclass[11pt,leqno]{amsart}

\usepackage{amsthm}
\usepackage{amsmath}
\usepackage{amssymb}
\usepackage[all,cmtip]{xy}
\usepackage[T1]{fontenc}
\usepackage{mathrsfs}
\usepackage{enumerate}
\theoremstyle{plain}
\newtheorem{thm}{Theorem}
\newtheorem{lem}{Lemma}
\newtheorem{prop}{Proposition}
\newtheorem{cor}{Corollary}
\theoremstyle{definition}

\newtheorem{rem}{Remark}

\title{Hilbert 90 for Algebras with Conjugation}
\date{}

\author{Patrik Lundström}

\address{Patrik Lundstr\"{o}m,
University West,
Department of Engineering Science,
Trollh\"{a}ttan,
Sweden}
\email{Patrik.Lundstrom@hv.se}


\begin{document}

\maketitle

\begin{abstract}
We show a version of Hilbert 90 that is valid for a large class of
algebras many of which are not commutative, distributive or associative.
This class contains the $n$th iteration
of the Conway-Smith doubling procedure.
We use our version of Hilbert 90 to parametrize all
solutions in ordered fields to the norm one equation for such algebras.
\end{abstract}

\section{Introduction}

Let $k$ be a field and suppose that $K$ is a finite Galois field
extension of $k$ with Galois group $G$. Recall that the induced norm
map $n : K \rightarrow k$ is defined by $n(a) = \prod_{\sigma \in G}
\sigma(a)$ for all $a \in K$. The well known Theorem 90 in Hilbert's
Zahlbericht \cite{hilbert}, also known as Hilbert 90, asserts that
if $G$ is cyclic, then an $a \in K$ satisfies $n(a)=1$ if and only
if for each automorphism $\sigma$ that generates $G$, there is a
nonzero $b \in K$ satisfying $\sigma(b)a=b$. For a proof, see e.g.
Chapter VI in \cite{lang}. The impetus for this article is the
observation that a version (see Theorem 1) of Hilbert 90 in degree
two holds for a large class of algebras which are not necessarily
commutative, distributive or associative. Namely, let $K$ be a
$k$-algebra. By this we mean that the following three properties should hold
(i) $K$ is a left $k$-vector space
equipped with a multiplication and a multiplicative identity 1;
(ii) $k$ is a subset of the center of $K$ and the additive and
multiplicative structure on $k$ equals the restriction of the
additive and multiplicative structure on $K$;
(iii) $a(bc)=(ab)c$, for $a,b,c \in K$,
whenever at least one of $a$, $b$ or $c$ belongs to $k$.
We say that $K$ is right (weak) distributive
if $(a+b)c = ac + bc$, for $a,b,c \in K$ ($a \in k$);  left (weak)
distributivity is analogously defined. Furthermore, we say that $K$
is (weak) distributive if it is both left and right (weak)
distributive.
Recall that $K$ is called left (or right) alternative if
$a(ab) = a^2 b$ (or $(ab)b=a b^2$), for $a,b \in K$;
$K$ is called alternative if it is both left
and right alternative.
Let ${}^{-} : K \rightarrow K$ be a self inverse
$k$-linear function which restricts to the identity map on $k$; we
will refer to such a map as a conjugation. The conjugation ${}^{-} :
K \rightarrow K$ is called an involution if it is a ring
antiautomorphism of $K$, that is, if $\overline{ab}=\overline{b}
\overline{a}$, for $a,b \in K$. The functions $n : K \rightarrow K$
and $t: K \rightarrow K$ defined by $n(a) = \overline{a}a$ and $t(a)
= a + \overline{a}$, for $a \in K$, respectively, will be referred
to as the norm and trace on $K$. We say
that an element $a \in K$ is imaginary if $t(a)=0$. We say that $n$
is anisotropic if $n(a) \neq 0$, for nonzero $a \in K$, and we say
that $n$ is multiplicative if $n(ab) = n(a)n(b)$, for $a,b \in K$.
Furthermore, we say that $n$ is symmetric if $n(\overline{a})=n(a)$,
for $a \in K$.

In Section \ref{algebraswithconjugation},
we show the following version of Hilbert 90
for algebras with conjugation.

\begin{thm}\label{hilbert90}
Suppose that $K$ is a weak distributive
left alternative $k$-algebra equipped with a conjugation
making the trace $k$-valued and the norm anisotro\-pic and $k$-valued.
If $K$ has a nonzero imaginary element,
then an $a \in K$ satisfies $n(a)=1$ if and
only if there is a nonzero $b \in K$ satisfying
$\overline{b}a = b$.
\end{thm}

There are lots of algebras with conjugation. The most well known
examples are the ones constructed via the Cayley-Dickson doubling
process (for more details, see e.g. \cite{bourbaki} or \cite{schafer94}).
Namely, suppose that $K$ is a $k$-algebra equipped with a conjugation and
that we have chosen a $D \in k$. Then there is a conjugation and a
multiplication on the $k$-vector space $K \times K$ defined by
the relations $\overline{(a,b)} = (\overline{a},-b)$
and $(a,b)(c,d) = (ac+Dd\overline{b} , \overline{a}d+cb)$,
respectively, for $a,b,c,d \in K$.
The resulting structure, called the Cayley-Dickson double of $K$ with
respect to $D$, is in this article denoted $K \times_D K$; it is a
$k$-algebra with a nonzero multiplicative identity $(1,0)$ and a
nonzero imaginary element $(0,1)$. Starting with the real numbers
equipped with the trivial conjugation this process successively,
with $D=-1$ at each step,
yields the complex numbers, the quaternions, the octonions and the
sedenions. However, their properties keep degrading; the complex
numbers are not ordered, the quaternions are not commutative, the
octonions are not associative and the sedenions contain zero
divisors and hence do not have a multiplicative norm.

The absence of a multiplicative norm for the sedenions
is implied by the well known Hurwitz's theorem \cite{hurwitz}
which states that any real distributive algebra
equipped with a multiplicative norm is
necessarily isomorphic to either the real numbers,
the complex numbers, the quaternions or
the octonions.
The existence of these four algebras is
equivalent to that of $N$-square identities
$(x_1^2 + \cdots + x_N^2)
(y_1^2 + \cdots + y_N^2) =
(z_1^2 + \cdots + z_N^2)$
for $N = 1,2,4,8$ in which the $z_k$ are bilinear functions of the
$x_i$ and $y_j$.
In 1967 Pfister \cite{pfister, pfister95}, quite surprisingly, showed
the existence of such $N$-square identities for any
$N$ which is a power of two provided that
the $z_k$ are rational functions of the $x_i$ and $y_j$.
In his book \cite{con03} Conway pointed out that
such identities can be produced by various
modifications of the Cayley-Dickson procedure.
J. D. H. Smith \cite{smith95} has proposed a similar doubling
procedure, which, however, does not behave well
algebraically; e.g. when applied to the quaternions
it does not yield the octonions.
Here, we will be concerned with
a variation of Conway's procedure proposed by W. D. Smith
in his extensive article \cite{smith04}.
Since it seems that it was developed, at least
partially, in collaboration with Conway I will
refer to this procedure as Conway-Smith doubling
(for more details concerning priority, see Section 32 in loc. cit.).
The Conway-Smith process starting from the real numbers, not only produces the
complex numbers, the quaternions and the octonions, but keeps going
to yield $2^n$-dimensional algebraical structures with
multiplicative norms for all positive integers $n$; Smith calls
these structures the $2^n$-ons. This
procedure supposes that $K$ is a division algebra over $k$. By this
we mean that to each nonzero $a \in K$ there is a unique $b \in K$
satisfying $ab = ba = 1$; in that case we will write $b = a^{-1}$.
The Conway-Smith double of $K$ with respect to
$D \in k$ (Smith loc. cit. always keeps $D=-1$),
in this article denoted $K \times^D K$, is the $k$-vector space
$K \times K$ equipped with the same conjugation as for the
Cayley-Dickson double but with a multiplication defined,
for all $a,b,c,d \in K$, by
$$(a,b)(c,d) =
\left( ac + D \overline{ b \overline{d} } , \overline{\overline{b}
\overline{c}} + \overline{ \overline{b} \overline{ \overline{a}
\overline{ \overline{b^{-1}} \ \overline{d}} } }    \right)$$ if $b
\neq 0$ and $$(a,b)(c,d)=(ac,\overline{a}d)$$ if $b=0$.
The $k$-algebra $K \times^D K$ has a multiplicative identity
$(1,0)$ and a nonzero imaginary element $(0,1)$.

In Section \ref{cayleydickson} and Section \ref{conwaysmith},
we investigate when Theorem 1 holds for
the Cayley-Dickson and Conway-Smith
doubling processes (see Corollaries 1, 2 and 3). As an
application, we parametrize explicitely the solutions to the norm
one equation for the $n$th iteration of the Conway-Smith doubling of
$k$ (see Remark 4). We are also able
to show that this works for the first, second and third
Cayley-Dickson doubling of $k$ (see Remark 3).
Along the way, we give an alternative proof of
norm multiplicativity of the $2^n$-ons (see Corollary 5)
and for the first, second and third Cayley-Dickson doubling (see Corollary 3).
This and a few other results of ours
can be found in Smith's article \cite{smith04} in the case when the
doubling is defined by $D=-1$. However, whereas Smith loc. cit. often resorts
to arguments depending on (quasi) matrix representations in his
proofs, we use algebraical properties of the doubling process alone.
For the convenience of the reader, we have tried to
make our presentation as self-contained as possible.

\section{Algebras with Conjugation}\label{algebraswithconjugation}

In the end of this section, we show Theorem \ref{hilbert90}.
Along the way, and for use in the following sections, we show
some results concerning general algebras with conjugation
(see Propositions \ref{weakhilbert90}-\ref{multiplicativenorm}).
For the rest of the article, we
suppose that $K$ is a $k$-algebra equipped
with a conjugation ${}^{-} : K \rightarrow K$.

\begin{prop}\label{weakhilbert90}
Suppose that $K$ is weak distributive
and equipped with a nonzero imaginary element.
If an $a \in K$ satisfies $n(a)=1$, then there is a nonzero
$b \in K$ satisfying $\overline{b}a=b$.
\end{prop}

\begin{proof}
Suppose that $n(a)=1$ for some $a \in K$.
We consider two cases.

Case 1: $a=-1$.
Let $b$ be a nonzero imaginary element of $K$. Then,
$\overline{b}=-b$ and hence,
by left weak distributivity, we get that
$\overline{b}a = (-b)(-1) = (-b)(-1) + 0 =
(-b)(-1) + (-b) + b =
(-b)(-1 + 1) + b = (-b) \cdot 0 + b = b$.

Case 2: $a \neq -1$. Let $b = a+1$.
Then, by right weak distributivity, we get that
$\overline{b}a = \left( \overline{a+1} \right) a =
\left( \overline{a} + 1 \right) a = \overline{a}a+a =
n(a)+a = 1+a = b$.
\end{proof}

\begin{prop}\label{leftalternativenorm}
If $K$ is weak right distributive and
the trace and norm on $K$ are $k$-valued, then $K$
is left alternative if and only if
$\overline{a}(ab) = n(a)b$ for all $a,b \in K$.
In that case, if the norm is anisotropic
and symmetric, then $K$ is a division algebra
with $a^{-1} = \overline{a}n(a)^{-1}$ for all
nonzero $a \in K$.
\end{prop}

\begin{proof}
Suppose that $K$ is weak right distributive and
that the trace and norm are $k$-valued.
Take $a,b \in K$.
First we show the ''only if'' part of the proof.
Suppose that $K$ is alternative. Then
$$\overline{a}(ab) = (t(a)-a)(ab) =
t(a)ab - a(ab) = t(a)ab - a^2 b =$$
$$= t(a)ab - ( (t(a)-\overline{a})a )b =
t(a)ab - ( t(a)a - \overline{a}a )b =$$
$$= t(a)ab - ( t(a)a - n(a) )b =
t(a)ab - (t(a)ab - n(a)b) = n(a)b.$$
Therefore, $\overline{a}(ab) = n(a)b$.
Now we show the ''if'' part of the proof.
Suppose that $\overline{a}(ab) = n(a)b$. Then
$$a(ab) = (t(a) - \overline{a})(ab) =
t(a)ab - \overline{a}(ab) = t(a)ab - n(a)b =$$
$$= (t(a)a - n(a))b = (t(a)a - \overline{a}a)b =
((t(a)-\overline{a})a)b = a^2b.$$
Therefore, $a(ab)=a^2b$.
Now suppose that the norm also is anisotropic
and symmetric and that $a$ is nonzero.
Then $$\overline{a}n(a)^{-1}a = \overline{a} a n(a)^{-1}
= n(a) n(a)^{-1} = 1$$
and $$a \overline{a} n(a)^{-1} =
n(\overline{a})n(a)^{-1} = n(a) n(a)^{-1} = 1.$$
Therefore, $\overline{a} n(a)^{-1}$ is a
multiplicative inverse of $a$.
Suppose that $ab=1$. Then
$$b = n(a) b n(a)^{-1} = \overline{a}(ab) n(a)^{-1}
= \overline{a} \cdot 1 \cdot n(a)^{-1} = \overline{a}n(a)^{-1}.$$
Therefore, the multiplicative inverse of $a$ is unique.
\end{proof}

\begin{prop}\label{respectssquares}
Suppose that the trace on $K$ is $k$-valued.
(a) If $K$ is weak right distributive,
then $t(a)b = ab + \overline{a}b$ for all $a,b \in K$;
(b) If $K$ is weak distributive, then
the norm on $K$ is symmetric;
(c) If $K$ is weak distributive and the norm on $K$
is $k$-valued, then the conjugation
on $K$ respects squares;
(d) If $K$ is left alternative and weak right distributive
and the norm is $k$-valued, then the norm respects squares.
\end{prop}

\begin{proof}
Take $a,b \in K$ and suppose that the trace
on $K$ is $k$-valued.

(a) If $K$ is weak right distributive, then
$$t(a)b = ab + t(a)b - ab = ab + (t(a)-a)b = ab + \overline{a}b.$$

(b) If $K$ is weak distributive, then
$$n(\overline{a}) = a \overline{a} = a(t(a)-a) = t(a)a - a^2
= (t(a)-a)a = \overline{a}a = n(a).$$

(c) If $K$ is weak distributive and the norm on $K$ is $k$-valued, then
$$\overline{a^2} = \overline{(t(a)-\overline{a}) a} =
\overline{ t(a)a - n(a) } = t(a)\overline{a} - n(a) =$$
$$= \overline{a}t(a) - \overline{a}a =
\overline{a}( t(a) - a ) = \overline{a}^2.$$

(d) By Proposition \ref{leftalternativenorm}, left alternativity
of $K$ and (c), we get that
$$n(a^2) = \overline{a^2} a^2 = \overline{a}^2 a^2 =
\overline{a} ( \overline{a} a^2 ) =
\overline{a} ( n(a) a) = n(a) \overline{a}a = n(a)^2.$$
\end{proof}

In many cases the norm not only respects squares,
but is in fact multiplicative. In the sequel,
we will use the following result.

\begin{prop}\label{multiplicativenorm}
Suppose that $K$ is right weak distributive and left alternative
with $k$-valued trace and norm.
If there is a $k$-linear map $T : K \rightarrow K$
which restricts to the identity map on $k$ and
satisfies
\begin{equation}\label{Tidentity}
T \left( \overline{a(bc)} \right) =
T \left( \overline{c} \left( \overline{b} \overline{a} \right) \right)
\end{equation}
for all $a,b,c \in K$, then the norm on $K$ is multiplicative.
\end{prop}

\begin{proof}
Take $a,b \in K$. Then, by (\ref{Tidentity}) and
Proposition \ref{leftalternativenorm}, we get that
$$n(ab) = T(n(ab)) =
T \left( \overline{ab} \ ab \right) =
T \left( \overline{ab} \left( \overline{\overline{a}} \overline{\overline{b}} \right) \right)=
T \left( \overline{ \overline{b} \left( \overline{a}(ab) \right) } \right) =$$
$$= T \left( \overline{ \overline{b} n(a) b } \right) =
n(a) T \left( \overline{ \overline{b}b } \right)
= n(a) T \left( \overline{n(b)} \right) = n(a)n(b).$$
\end{proof}

\begin{rem}\label{classicalargument}
If $K$ is distributive and alternative
and equipped with an involution making
the trace and the norm $k$-valued, then,
by a classical argument,
the norm on $K$ is multiplicative.
Namely, by a theorem of
Artin (see Theorem 3.1 in \cite{schafer94}),
every subalgebra generated by two elements
of a distributive alternative algebra is associative.
Therefore, for any $a,b \in K$, we get that
$$n(ab) =
\overline{(ab)}(ab) =
(\overline{b}\overline{a})(ab) =
(t(b)-b)(t(a)-a)ab=$$
$$=(t(b)-b)( (t(a)-a)a )b =
\overline{b}(\overline{a}a)b=
\overline{b} n(a) b =
n(a) \overline{b} b =
n(a)n(b).$$
Since many of the structures we study in
the sequel possess conjugations which are not in general
involutions, Proposition \ref{multiplicativenorm} is still motivated.
\end{rem}

\subsection*{Proof of Theorem \ref{hilbert90}.}
The ''only if'' part of the claim follows from Proposition 1.
Now we show the ''if'' part of the claim. Suppose that
$\overline{b}a = b$ for some nonzero $a,b \in K$.
Since $K$ is left alternative and the norm is symmetric,
we get, by Proposition \ref{leftalternativenorm}, that
$n(b)a = n(\overline{b})a = b(\overline{b}a)= b^2$.
By taking norms we get, by Proposition \ref{respectssquares}(d),
that $n(b)^2 n(a) = n(b)^2$. Since $n(b)$ is a nonzero
element of $k$, this implies that $n(a) = 1$. {\hfill $\square$}

\section{Cayley-Dickson doubling}\label{cayleydickson}

In this section, we investigate when Theorem \ref{hilbert90} holds for algebras
gene\-rated by the Cayley-Dickson doubling process (see Corollaries
\ref{cayleydicksonweakhilbert90} and \ref{cayleydicksonhilbert90}).
For the rest of the article, we assume that $D$ is an
element of $k$.

\begin{prop}\label{cayleydicksondistributive}
(a) $K$ is weak right distributive if and only if $K \times_D K$ is weak right
distributive; (b) $K \times_D K$ is (weak) left distributive if and
only if $K \times_D K$, and hence $K$, is (weak) distributive;
(c) $K \times_D K$ is right distributive if and
only if $K \times_D K$, and hence $K$, is distributive;
(d) $K$ is weak distributive if and only if $K \times_D K$ is weak distributive.
\end{prop}

\begin{proof}
(a) The ''if'' part of the claim is clear since $K$ is contained in
$K \times_D K$. Now we show the ''only if'' part of the claim.
Suppose that $K$ is weak right distributive. Take $a_1 \in k$ and
$a_2,b,c,d \in K$. Then
$$(a_1+a_2,b)(c,d) = ( (a_1+a_2)c + Dd\overline{b} , (\overline{a_1+a_2})d + cb ) =$$
$$=( (a_1+a_2)c + Dd\overline{b} , (\overline{a}_1 + \overline{a}_2)d + cb ) =
( a_1 c + a_2 c + Dd\overline{b} , \overline{a}_1 d + \overline{a}_2
d + cb ) =$$
$$= (a_1 c , \overline{a}_1 d) + (a_2 c + Dd\overline{b}
, \overline{a}_2 d + cb ) = (a_1,0)(c,d) + (a_2,b)(c,d).$$ Therefore
$K \times_D K$ is weak right distributive.

(b) Suppose that $K \times_D K$ is left distributive. Since $K$ is
contained in $K \times_D K$, we get that $K$ is left distributive.
Now we show that $K$ is right distributive. Take $a,b,c \in K$. Then
$$(0,0) = (0,c)(a+b,0) - (0,c)(a,0)-(0,c)(b,0) = (0,(a+b)c - ac -bc).$$
Therefore $0 = (a+b)c - ac - bc$. Now we show that $K \times_D K$,
and hence $K$, is right distributive. Take $a,b,c,d,e,f \in K$. Then
$$( (a,b) + (c,d) ) (e,f) = (a+c,b+d)(e,f) =$$
$$=( (a+c)e + Df(\overline{b+d}) , (\overline{a+c})f + e(b+d) ) $$
$$=( ae+ce + Df(\overline{b}+\overline{d}) , (\overline{a}+\overline{c})f + eb+ed ) $$
$$=( ae+ce + D f \overline{b}+ D f \overline{d} , \overline{a}f + \overline{c}f + eb+ed ) $$
$$=( ae + D f \overline{b} , \overline{a}f +  eb) +
( ce + D f \overline{d} , \overline{c}f + ed )
= (a,b)(e,f) + (c,d)(e,f).$$
The weak part of (b) is proved in a similar way.

(c) This can be proved in a fashion analogous to (b).

(d) This follows from (a) and (b).
\end{proof}

\begin{cor}\label{cayleydicksonweakhilbert90}
Suppose that $k$ is equipped with the identity conjugation. If we
for each positive integer $i$ choose $D_i \in k$ and recursively
define the $k$-algebras $k_i$, for $i \geq 0$, by $k_0 = k$ and
$k_i = k_{i-1} \times_{D_i} k_{i-1}$, for $i \geq 1$, then, for all $i
\geq 1$, if an $a \in k_i$ satisfies $n(a)=1$, then there is a
nonzero $b \in k_i$ satisfying $\overline{b}a=b$.
\end{cor}

\begin{proof}
This follows immediately from Propositions \ref{weakhilbert90}
and \ref{cayleydicksondistributive}.
\end{proof}

\begin{prop}\label{cayleydicksoninvolution}
(a) The conjugation on $K$ is an involution if and only if the
conjugation on $K \times_D K$ is an involution; (b) The trace on $K$
is $k$-valued if and only if the trace on $K \times_D K$ is
$k$-valued; (c) The norm on $K$ is symmetric if and only if the norm
on $K \times_D K$ is symmetric. In that case (i) the norm on $K$ is
$k$-valued if and only if the norm on $K \times_D K$ is $k$-valued;
(ii) the norm on $K \times_D K$ is anisotropic if and only if the
norm on $K$ is anisotropic and $D$ is an element of $k$ not
belonging to the set of quotients $n(a)/n(b)$, for $a \in K$ and
nonzero $b \in K$.
\end{prop}

\begin{proof}
All the ''only if'' parts of the proof follow
from the fact that $K$ is contained in $K \times_D K$.
Now we show the ''if'' parts of the proof.

(a) Suppose that the conjugation on $K$
is an involution. Take $a,b,c,d \in K$. Then
$$\overline{(a,b)(c,d)}=
\overline{ \left( ac+Dd\overline{b} , \overline{a}d+cb \right) } =
\left( \overline{ac}+D \overline{d \overline{b}} ,
-\overline{a}d-cb \right) =$$
$$= \left( \overline{c} \ \overline{a} + D
\overline{\overline{b}} \ \overline{d}, -\overline{a}d-cb \right)=
\left( \overline{c} \ \overline{a} + D(-b)(\overline{-d}),
\overline{\overline{c}}(-b) + \overline{a}(-d) \right) =$$
$$= (\overline{c},-d) (\overline{a},-b) = \overline{(c,d)} \
\overline{(a,b)}.$$
Therefore the conjugation on $K \times_D K$
is an involution.

(b) Suppose that the trace on $K$ is $k$-valued.
Take $a,b \in K$. Then
$$t(a,b) = (a,b) + \overline{(a,b)} =
(a,b) + (\overline{a},-b) = (a + \overline{a},b-b)=
(t(a),0)$$ which belongs to $k$.

(c) Suppose that the norm on $K$ is symmetric.
Take $a,b \in K$. Then
$$n((a,b)) = \overline{(a,b)} (a,b) =
(\overline{a},-b)(a,b) =
(\overline{a}a + Db(-\overline{b}) ,
\overline{\overline{a}}b + a(-b) ) =$$
$$= (n(a) + D(-b)\overline{b} , ab - ab)
= ( n(\overline{a}) + D(-b)\overline{b} , 0) =$$
$$= ( \overline{\overline{a}} \ \overline{a} + D(-b)\overline{b} ,
\overline{a}(-b) + \overline{a}b) =
( a\overline{a} + D(-b)\overline{b} ,
\overline{a}(-b) + \overline{a}b) ) =$$
$$= (a,b)(\overline{a},-b) =
(a,b) \overline{(a,b)} =
\overline{\overline{(a,b)}} (a,b) =
n \left( \overline{(a,b)} \right).$$
Therefore, the norm on $K \times_D K$ is symmetric.
Furthermore, the above calculation shows that
$n((a,b)) = n(a) - Dn(b)$, for $a,b \in K$.
From this equality, (i) and (ii) follow immediately.
\end{proof}

\begin{prop}\label{cayleydicksonalternative}
If $K$ is left distributive,
weak right distributive and both
the trace and norm on $K$ are $k$-valued,
then $K$ is associative if and only if
$K \times_D K$ is left alternative.
\end{prop}

\begin{proof}
Suppose that $K \times_D K$, and hence $K$, is left alternative.
Take $a,b,c \in K$. Then
$$(0,0) = (a,c)^2 (b,0) - (a,c)( (a,c)(b,0) ) =$$
$$=(a^2 + D c \overline{c} , \overline{a}c+ac)(b,0) - (a,c)(ab,bc)$$
$$=(a^2 + Dn(c) , t(a)c)(b,0) - (a^2 b + D(bc)\overline{c} , \overline{a}bc + (ab)c)$$
$$= (a^2 b + Dn(c)b , t(a)bc) - (a^2 b + Dbn(c) , \overline{a}bc + (ab)c )$$
$$= (0 , (t(a) - \overline{a})bc - (ab)c ) = (0 , a(bc) - (ab)c).$$
Therefore $a(bc)-(ab)c=0$ and hence $K$ is associative.

Suppose that $K$ is associative.
Take $a,b,c,d \in K$. Then, by Proposition \ref{respectssquares}(c),
we get that
$$ (a,b)^2 (c,d) =  (a^2+Dn(b) , t(a)b)(c,d) =$$
$$= ( (a^2+Dn(b))c + Dt(a)d\overline{b} ,
(\overline{a^2+Dn(b)})d + t(a)cb )$$
$$= (a^2 c + Dn(b)c + Dt(a)d\overline{b} ,
\overline{a}^2 d + Dn(b)d + t(a)cb )$$
and, by Proposition \ref{respectssquares}(a), we get that
$$ (a,b)((a,b)(c,d) ) = (a,b) (ac+Dd\overline{b} , \overline{a}d+cb) =$$
$$= ( a(ac+Dd\overline{b}) +
D(\overline{a}d+cb)\overline{b} ,
\overline{a}(\overline{a}d+cb)+(ac+Dd\overline{b})b )$$
$$= (a^2 c + Dad\overline{b} + D\overline{a}d\overline{b} + Dcn(b) ,
\overline{a}^2 d + \overline{a}cb + acb + Ddn(b) )$$
$$= (a^2 c + Dt(a)d\overline{b} + Dcn(b) ,
\overline{a}^2 d + Ddn(b) + t(a)cb).$$
Therefore $(a,b)^2 (c,d) = (a,b)((a,b)(c,d) )$
and hence $K \times_D K$ is alternative.
\end{proof}

Recall that $k$ is called ordered if it is equipped
with a total order $\leq$ satisfying the following
two properties for all $a,b,c \in k$:
(i) if $a \leq b$, then $a+c \leq b+c$;
(ii) if $0 \leq a$ and $0 \leq b$, then $0 \leq ab$.
An element $a \in k$ is called positive if
$0 \leq a$ and $0 \neq a$. %; this is denoted $0 < a$.
Negative elements of $k$ are analogously defined.

\begin{cor}\label{cayleydicksonhilbert90}
Suppose that the $k$-algebras
$k_1$, $k_2$ and $k_3$ are defined as in Corollary 1.
Then an $a \in k_i$ satisfies $n(a)=1$
if and only if there is a nonzero $b \in k_i$
with $\overline{b}a=b$ if
(a) $i=1$ and $D_1$ is not a square in $k$ or
(b) $i=1,2,3$ whenever $k$ is an ordered field and
$D_1$, $D_2$ and $D_3$ are negative.
\end{cor}

\begin{proof}
This follows immediately from Theorem \ref{hilbert90} and
Propositions \ref{cayleydicksondistributive}-\ref{cayleydicksonalternative}.
\end{proof}

\begin{rem}
It is not clear to the author at present
whether the conclusion of Corollary \ref{cayleydicksonhilbert90} holds
for $i \geq 4$ since the techniques we have used
depend on the fact that the algebras
are left alternative; it is well known that
this is not satisfied for algebras containing
the sedenions, that is, if $i \geq 4$.
\end{rem}

\begin{rem}\label{cayleydicksonnormone}
By folklore, Hilbert 90 can be used to parametrize solutions to the
norm one equation. Now we apply this idea to the $k$-algebras $k_1$,
$k_2$ and $k_3$ as defined in Corollary \ref{cayleydicksonweakhilbert90}.
Put $C_i = -D_i$ for
$i = 1,2,3$. Given $x \in k_i$ and $1 \leq i \leq 3$, then, by the
proof of Theorem \ref{hilbert90},
the equality $n(x)=1$ holds if and only if there is a
nonzero $s \in k_i$ with
\begin{equation}\label{decomposition}
x = \frac{s^2}{n(s)} = \frac{2s_1^2-n(s)}{n(s)} + \frac{2s_1(s-s_1)}{n(s)}
\end{equation}
where $s_1$ is the $k$-part of $s$. Thus, by calculating the norm of
a general element of $k_i$, for $i = 1,2,3$, we can parametrize the
solutions to the norm one equation.
We now make this parametrization explicit in
the cases $i = 1,2,3$ separately.

If $-C_1$ is not a square in $k$, then, by
(\ref{decomposition}),  we get that $x_1,x_2 \in
k$ satisfy
$$x_1^2 + C_1 x_2^2 = 1$$ if and only if there are $s_1,s_2 \in k$, not
both zero, such that
$$x_1 = \frac{ s_1^2 - C_1 s_2^2 }{ s_1^2 + C_1 s_2^2 } \quad \quad
\mbox{and} \quad \quad x_2 = \frac{2 s_1 s_2}{ s_1^2 + C_1 s_2^2 }$$
In particular, this holds if $k$ is an ordered field
and $C_1$ is positive.

If $k$ is an ordered field and $C_1$ and $C_2$ are positive,
then, by (\ref{decomposition}),
we get that $x_i \in k$, for $i=1,2,3,4$, satisfy
$$x_1^2 + C_1 x_2^2 + C_2 x_3^2 + C_1 C_2 x_4^2 = 1$$ if and
only if there are $s_1,s_2,s_3,s_4 \in k$, not all zero, such that
$$x_1 = \frac{s_1^2 - C_1 s_2^2 - C_2 s_3^2 - C_1 C_2 s_4^2}{s_1^2 + C_1 s_2^2 + C_2 s_3^2 + C_1 C_2
s_4^2}$$ and
$$x_i = \frac{2s_1 s_i}{s_1^2 + C_1 s_2^2 + C_2 s_3^2 + C_1 C_2
s_4^2}$$ for $i = 2,3,4$.

If $k$ is an ordered field and $C_1$, $C_2$ and
$C_3$ are positive, then, by (\ref{decomposition}),
we get that $x_i \in k$, for $i=1,\ldots,8$,
satisfy {\small
$$x_1^2 + C_1 x_2^2 + C_2 x_3^2 + C_1 C_2 x_4^2 +
C_3 x_5^2 + C_1 C_3 x_6^2 + C_2 C_3 x_7^2 + C_1 C_2 C_3 x_8^2 = 1$$}
if and only if there are $s_i \in k$, for $1 \leq i \leq 8$, not all
zero, such that {\small
$$x_1 = \frac{s_1^2 - C_1 s_2^2 - C_2 s_3^2 - C_1 C_2 s_4^2 -
C_3 s_5^2 - C_1 C_3 s_6^2 - C_2 C_3 s_7^2 - C_1 C_2 C_3 s_8^2}{s_1^2
+ C_1 s_2^2 + C_2 s_3^2 + C_1 C_2 s_4^2 + C_3 s_5^2 + C_1 C_3 s_6^2
+ C_1 C_3 s_7^2 + C_1 C_2 C_3 s_8^2}$$ } and {\small
$$x_i = \frac{2 s_1 s_i}{s_1^2 + C_1 s_2^2 + C_2 s_3^2 + C_1 C_2 s_4^2 + C_3 s_5^2 + C_1
C_3 s_6^2 + C_2 C_3 s_7^2 + C_1 C_1 C_3 s_8^2}$$ } for $2 \leq i
\leq 8$.
\end{rem}

We end this section with a digression on norm
multiplicativity of Cayley-Dickson doubles.

\begin{prop}\label{cayleydicksonT}
Suppose that $K$ is associative and left distributive
and there is a $k$-linear map
$T : K \rightarrow K$ satisfying the following three conditions
(i) $T \left (\overline{a} \right) = T(a)$;
(ii) $T(ab) = T(ba)$;
(iii) $T \left( \overline{a(bc)} \right) =
T \left( \overline{c} \left( \overline{b}\overline{a} \right) \right)$
for all $a,b,c \in K$.
If we extend the map $T$ to $K \times_D K$
by the relation $T((a,b)) = T(a)$, for $a,b \in K$,
then (i), (ii) and (iii) hold for all
$a,b,c \in K \times_D K$.
In that case, if $K$ is
weak right distributive, then
the norm on $K \times_D K$ is multiplicative.
\end{prop}

\begin{proof}
Take $a = (a_1,a_2)$, $b = (b_1,b_2)$ and $c = (c_1,c_2)$
in $K \times_D K$.

First we show (i):
$T(\overline{a}) = T(\overline{a}_1,-a_2) =
T(\overline{a}_1) = T(a_1) = T(a)$.

Next we show (ii):
$T(ab) = T(a_1 b_1 + D b_2 \overline{a}_2) =
T(a_1 b_1) + D T(b_2 \overline{a}_2) =
T(b_1 a_1) + D T \left( \overline{b_2 \overline{a}_2} \right)
= T(b_1 a_1) + D T(a_2 \overline{b}_2) = T(ba)$.

Finally, we show (iii). First note that left distributivity in
combination with (ii) implies that
$$T((x+y)z) = T(z(x+y)) = T(zx+zy) =$$
$$= T(zx) + T(zy) =  T(xy) + T(yz) = T(xy+yz)$$
for all $x,y,z \in K$.
Therefore, by (i), we get that
$$T \left( \overline{a(bc)} \right) =
T \left(  a_1 \left( b_1 c_1 + D c_2 \overline{b}_2 \right) \right) +
DT \left( \left( \overline{b}_1 c_2+c_1 b_2 \right) \overline{a}_2  \right) =$$
$$= T \left( \overline{a_1 (b_1 c_1)} \right) + DT \left( \overline{ a_1 \left( c_2 \overline{b}_2 \right)}\right) +
DT \left( \overline{ (\overline{b}_1 c_2) \overline{a}_2 } \right) +
DT \left( \overline{ (c_1 b_2) \overline{a}_2 } \right)  $$
$$= T \left( \overline{c}_1 \left( \overline{b}_1 \overline{a}_1 \right) \right) +
DT \left( b_2 \left( \overline{c}_2 \ \overline{a}_1 \right) \right) +
DT \left( a_2 \left( \overline{c}_2 b_1 \right) \right) +
DT \left( a_2 \left( \overline{b}_2 \overline{c}_1 \right) \right).$$
Likewise
$$T \left( \overline{c} \left( \overline{b} \overline{a} \right) \right) =
T \left( \overline{c}_1 \left( \overline{b}_1 \overline{a}_1 \right) \right) +
DT \left( \overline{c}_1 \left( a_2 \overline{b}_2 \right) \right) +
DT \left( (b_1 a_2) \overline{c}_2 \right) +
DT \left( \left( \overline{a}_1 b_2 \right) \overline{c}_2 \right)$$
which, by (ii) and associativity of $K$,
equals $T \left( \overline{a(bc)} \right)$.

The last part follows from the above and
Propositions \ref{multiplicativenorm} and \ref{cayleydicksonalternative}.
\end{proof}

\begin{cor}\label{cayleydicksonmultiplicativenorm}
The norms on the $k$-algebras $k_1$, $k_2$ and $k_3$,
as defined in Corollary \ref{cayleydicksonweakhilbert90},
are multiplicative.
\end{cor}

\begin{proof}
Suppose that we define the $k$-linear map $T_3 : k_3 \rightarrow k$,
as the projection on $k$. By restriction, this
induces $k$-linear maps $T_2 : k_2 \rightarrow k$,
$T_1 : k_1 \rightarrow k$ and $T_0 : k_0 \rightarrow k$.
Since obviously $T_0 = {\rm id}_k$ satisfies (i), (ii)
and (iii) from Proposition \ref{cayleydicksonT}, it follows
that the same is true for $T_1$, $T_2$ and $T_3$.
Therefore, by Propositions 4, 5 and 7,
the norm is multiplicative on $k_1$, $k_2$ and $k_3$.
\end{proof}

\section{Conway-Smith doubling}\label{conwaysmith}

In this section, we investigate when Theorem 1 holds for algebras
gene\-rated by the Conway-Smith doubling process.

\begin{prop}\label{conwaysmithdoubling}
If $K$ is left alternative, left distributive,
right weak distributive, and has $k$-valued
trace and anisotropic $k$-valued norm,
then $K \times^D K$ is left
alternative, left distributive,
right weak distributive and has $k$-valued
trace and norm.
In that case, if $D$ does not belong to the set of quotients
$n(a)/n(b)$, for $a \in K$ and nonzero $b \in K$,
then $K \times^D K$ has anisotropic norm.
\end{prop}

\begin{proof}
Take $a,b,c,d,e,f \in K$.
First we show that $K \times^D K$ is left distributive.
Case 1: $b=0$. Then
$$ (a,b) ( (c,d)+(e,f) ) = (a,b)(c+e,d+f) = (a,0)(c+e,d+f)= $$
$$= ( a(c+e) , \overline{a}(d+f) ) = ( ac+ae , \overline{a}d + \overline{a}f )=
(ac,\overline{a}d) + (ae,\overline{a}f) =$$
$$=(a,0)(c,d) + (a,0)(e,f) =
(a,b)(c,d) + (a,b)(e,f).$$
Case 2: $b \neq 0$. Then
$$ (a,b) ( (c,d)+(e,f) ) = (a,b)(c+e,d+f) = $$
$$= \left( a(c+e) + D \overline{ b (\overline{d+f}) } ,
\overline{ \overline{b} (\overline{c+e}) } +
\overline{ \overline{b} \overline{ \overline{a}
\overline{ \overline{b^{-1}}  (\overline{d+f})} } } \right) $$
$$= \left( ac + D \overline{ b \overline{d} } ,
\overline{ \overline{b} \overline{c} } +
\overline{ \overline{b} \overline{ \overline{a}
\overline{ \overline{b^{-1}}  \overline{d}} } } \right) +
    \left( ae + D \overline{ b \overline{f} } ,
\overline{ \overline{b} \overline{e} } +
\overline{ \overline{b} \overline{ \overline{a}
\overline{ \overline{b^{-1}}  \overline{f}} } } \right) $$
$$= (a,b)(c,d) + (a,b)(e,f).$$
Now we show that $K \times^D K$ is weak right distributive.
If $a \in k$, then
$$( (a,0) + (c,d) )(e,f) = (a+c,d)(e,f) =$$
$$= \left( (a+c)e + D \overline{ d \overline{f} } ,
\overline{ \overline{d} \overline{e} } +
\overline{ \overline{d} \overline{ (\overline{a+c})
\overline{ \overline{d^{-1}}  \overline{f} } } } \right) $$
$$= \left( ae + ce + D \overline{ d \overline{f} } ,
\overline{ \overline{d} \overline{e} } +
\overline{ \overline{d} \overline{ (\overline{a}+\overline{c})
\overline{ \overline{d^{-1}}  \overline{f} } } } \right) $$
$$= \left( ae + ce + D \overline{ d \overline{f} } ,
\overline{ \overline{d} \overline{e} } +
\overline{
\overline{d}
\left(
\overline{
 \overline{a} \overline{ \overline{d^{-1}}  \overline{f} } } +
 \overline{
 \overline{c} \overline{ \overline{d^{-1}}  \overline{f} } }
\right)
}
\right) $$
$$= \left( ae + ce + D \overline{ d \overline{f} } ,
\overline{ \overline{d} \overline{e} }
+
\overline{
\overline{d}
\overline{
 \overline{a} \overline{ \overline{d^{-1}}  \overline{f} } }
}
+
\overline{
\overline{d}
 \overline{
 \overline{c} \overline{ \overline{d^{-1}}  \overline{f} } }
}
\right) $$
$$= \left( ae + ce + D \overline{ d \overline{f} } ,
\overline{ \overline{d} \overline{e} }
+
\overline{a}
\overline{
\overline{d} \
\overline{
 d^{-1}} \ \overline{f}
}
+
\overline{
\overline{d}
 \overline{
 \overline{c} \overline{ \overline{d^{-1}}  \overline{f} } }
}
\right) $$
$$= \left( ae + ce + D \overline{ d \overline{f} } ,
\overline{ \overline{d} \overline{e} }
+
\overline{a} f
+
\overline{
\overline{d}
 \overline{
 \overline{c} \overline{ \overline{d^{-1}}  \overline{f} } }
}
\right) $$
$$= (ae,\overline{a}f) +
\left( ce + D \overline{ d \overline{f} } ,
\overline{ \overline{d} \overline{e} }
+
\overline{
\overline{d}
 \overline{
 \overline{c} \overline{ \overline{d^{-1}}  \overline{f} } }
}
\right)
=(a,0)(e,f) + (c,d)(e,f).$$
Now we show that $K \times^D K$ is left alternative.

\vspace{1mm}

\noindent Case 1: $b=0$. Then, by Proposition \ref{respectssquares}(c), we get that
$$(a,b)^2 (c,d) = (a,0)^2(c,d) = (a^2,0) (c,d) =
(a^2 c , \overline{a^2}d) = (a^2 c , \overline{a}^2d) =$$
$$= (a(ac) , \overline{a}(\overline{a}d)) =
(a,0)(ac,\overline{a}d) = (a,0)( (a,0)(c,d) ) =
(a,b)( (a,b)(c,d) )$$
Case 2: $b \neq 0$ and $t(a) \neq 0$. Then, by Proposition \ref{respectssquares}(c),
we get that
$$(a,b)^2 (c,d) =
(a^2 + Dn(b) , \overline{\overline{b}\overline{a}} +
\overline{ \overline{b}\overline{\overline{a}\overline{\overline{b^{-1}}\overline{b}  } } })
(c,d) =$$
$$= (a^2+Dn(b) , \overline{\overline{b}\overline{a}} +
\overline{\overline{b}a})(c,d) =
(a^2 + Dn(b) , \overline{ \overline{b}\overline{a} + \overline{b}a })(c,d) =$$
$$= (a^2 + Dn(b) , \overline{ \overline{b}(\overline{a} + a) })(c,d) =
(a^2 + Dn(b) , t(a)b)(c,d) =$$
$$= ( (a^2+Dn(b))c + D \overline{bt(a)\overline{d}} ,
\overline{\overline{b}t(a)\overline{c}} +
\overline{ \overline{b}t(a) \overline{
(\overline{a^2+Dn(b)}) \overline{ \overline{b^{-1}t(a)^{-1}}\overline{d} } } } )$$
$$= ( a^2 c + Dn(b)c + Dt(a) \overline{b \overline{d}} ,
t(a)\overline{\overline{b}\overline{c}} +
\overline{ \overline{b} ( \overline{
\overline{a}^2 \overline{\overline{b^{-1}}\overline{d}} +
Dn(b) \overline{\overline{b^{-1}}\overline{d}} } ) } )$$
$$= ( a^2 c + Dn(b)c + Dt(a) \overline{b \overline{d}} ,
t(a)\overline{\overline{b}\overline{c}} +
\overline{ \overline{b} \overline{ \overline{a}^2
\overline{\overline{b^{-1}}\overline{d}} } } + Dn(b)d )$$
and
$$(a,b) ( (a,b)(c,d) ) = (a,b) ( ac+D\overline{b\overline{d}} ,
\overline{\overline{b}\overline{c}} +
\overline{ \overline{b}
\overline{ \overline{a} \overline{ \overline{b^{-1}}\overline{d} } } } ) =$$
$$= \left( a(ac+D \overline{b\overline{d}} +
D\overline{ b(\overline{b}\overline{c}+
\overline{b}\overline{ \overline{a} \overline{ \overline{b^{-1}}\overline{d} } }) } ,
\overline{ \overline{b}( \overline{ac} + Db\overline{d} ) } +
\overline{ \overline{b} \overline{ \overline{a} \overline{ \overline{b^{-1}}
( \overline{b}\overline{c} +
\overline{b} \overline{ \overline{a} \overline{ \overline{b^{-1}}\overline{d} } } ) } } }   \right)$$
$$= \left( a^2 c + D a \overline{ b\overline{d} } +
Dn(b)c + Dn(b)\overline{a}\overline{ \overline{b^{-1}}\overline{d} } ,
\overline{ \overline{b}\overline{ac} } +
Dn(b)d + \overline{ \overline{b} \overline{ \overline{a}c } } +
\overline{ \overline{b} \overline{ \overline{a}^2 \overline{ \overline{b^{-1}}\overline{d} } } }  \right) $$
$$= \left( a^2 c + Dn(b)c + Dt(a) \overline{b \overline{d}} ,
t(a)\overline{\overline{b}\overline{c}} +
\overline{ \overline{b} \overline{ \overline{a}^2
\overline{\overline{b^{-1}}\overline{d}} } } + Dn(b)d \right).$$
Case 3: $b \neq 0$ and $t(a)=0$. Then, by Proposition \ref{respectssquares}(c)
and the beginning of the
calculation in Case 2, we get that
$$(a,b)^2 (c,d) = (a^2 + Dn(b) , 0)(c,d)
= \left( (a^2+Dn(b))c , \overline{(a^2+Dn(b))}d   \right) =$$
$$= \left( a^2 c + Dn(b)c , \overline{a}^2 d + Dn(b)d \right).$$
By the second part of the calculation in Case 2 and
the fact that $\overline{a}^2 = -n(a) \in k$, we get that
$$(a,b)((a,b)(c,d)) =
\left( a^2 c + Dn(b)c ,
\overline{ \overline{b} \overline{ \overline{a}^2
\overline{\overline{b^{-1}}\overline{d}} } } + Dn(b)d \right) =$$
$$= \left( a^2 c + Dn(b)c , \overline{a}^2 d + Dn(b)d \right).$$
The last part follows from the fact that
$\overline{(a,b)}(a,b) = (n(a)-Dn(b) , 0)$.
\end{proof}

\begin{cor}\label{conwaysmithhilbert90}
Suppose that $k$ is an ordered field equipped with the identity
conjugation. If we for each nonnegative integer $i$ choose a negative
$D_i \in k$ and recursively define the $k$-algebras $k^i$, for $i
\geq 0$, by $k^0 = k$ and $k^i = k^{i-1} \times^{D_i} k^{i-1}$, for $i
\geq 1$, then, for all $i \geq 1$, an $a \in k^i$ satisfies $n(a)=1$
if and only if there is a nonzero $b \in k^i$ satisfying
$\overline{b}a=b$.
\end{cor}

\begin{proof}
This follows immediately from Theorem \ref{hilbert90}
and Proposition \ref{conwaysmithdoubling}.
\end{proof}

\begin{rem}\label{conwaysmithnormone}
Now we use our version of Hilbert 90 to parametrize the solutions to
norm one equations for $2^n$-ons. Suppose that $k$ is an ordered field
and that the $k$-algebras
$k^i$, for $i \geq 0$, are defined as in Corollary \ref{conwaysmithhilbert90}.
Fix a nonnegative integer $n$. Now we define $2^n$ coordinates on $k^n$
recursively in the following way. Suppose that $x = (y,z) \in k^n =
k^{n-1} \times^{D_n} k^{n-1}$ and $1 \leq j \leq 2^n$. The $j$th
coordinate of $x$ is defined as the $j$th coordinate of $y$, if $1
\leq j \leq  2^{n-1}$, and as the $j$th coordinate of $z$, if
$2^{n-1} < j \leq 2^n$. Label the $2^n$ coordinates of $x$ as
$x_1,x_2,\ldots,x_{2^n}$. Put $C_i = -D_i$, for $i \geq 0$.
Define a total order $<$ on the $2^n$ elements of $\{ 0,1 \}^n$
by letting $(a_1,\ldots,a_n) < (b_1,\ldots,b_n)$ if there is
a positive integer $j$ with $a_j = 0$, $b_j=1$ and
$a_i = b_i$ for $i > j$.
Define $e(1),e(2),\ldots,e(2^n) \in \{ 0,1 \}^n$ by
$e(1) < e(2) < \cdots < e(2^n)$.
For each $i$, let $C^{e(i)}$ denote the product
$\prod_{j=1}^n C^{e(i)_j}_j$, where
$e(i) = (e(i)_1,e(i)_2,\ldots,e(i)_n)$.
Then, by equation (\ref{decomposition}), we get that
$$n(x) = \sum_{i=1}^{2^n} C^{e(i)} x_i^2 = 1$$
if and only if there are $s_j \in k$, for $1 \leq j \leq 2^n$, not all
zero, such that
$$x_1 = \frac{s_1^2 - \sum_{j=2}^{2^n} C^{e(j)} s_j^2}
{\sum_{j=1}^{2^n} C^{e(j)} s_j^2}$$
and
$$x_i = \frac{2 s_1 s_i}{\sum_{j=1}^{2^n} C^{e(j)} s_j^2}$$
for $2 \leq i \leq 2^n$.

A special case of the above discussion,
interesting on it's own right, is the following.
Elements $x_1,x_2,\ldots,x_{2^n}$ of an ordered
field $k$ satisfy
$$x_1^2 + x_2^2 + \cdots x_{2^n}^2 = 1$$
if and only if there are $s_j \in k$, for $1 \leq j \leq 2^n$, not all
zero, such that
$$x_1 = \frac{s_1^2 - (s_2^2 + s_3^2 + \cdots s_{2^n}^2)}
{s_1^2 + s_2^2 + \cdots s_{2^n}^2}$$
and
$$x_i = \frac{2 s_1 s_i}{s_1^2 + s_2^2 + \cdots s_{2^n}^2}$$
for $2 \leq i \leq 2^n$.
In particular, this implies that a vector
$$x=(x_1,x_2,\ldots,x_{2^n+1})$$ of $2^n+1$ integers
satisfy
$$x_1^2 + x_2^2 + \cdots + x_{2^n}^2 = x_{2^n+1}^2$$
if and only if $x$
is a rational multiple of a vector of the form
$$\big( s_1^2 - (s_2^2 + s_3^2 + \cdots s_{2^n}^2) , 2 s_1 s_2,
2 s_1 s_3 , \ldots , 2 s_1 s_{2^n} , s_1^2 + s_2^2 + \cdots + s_{2^n}^2 \big)$$
for some integers $s_1,s_2,\ldots,s_{2^n}$.
Note that this result is stated without proof on p. 72 in 
\cite{andres}.
\end{rem}

We end this section with a digression on
norm multiplicativity of Conway-Smith doubles.

\begin{lem}\label{normlemma}
If $T : K \rightarrow K$
is a $k$-linear map satisfying the following
three conditions
(i) $T(ab)=T(ba)$;
(ii) $T \left( \overline{a} \right) = T(a)$;
(iii) $T \left( \overline{a(bc)} \right) =
T \left( \overline{c} \left( \overline{b} \overline{a} \right) \right)$,
for all $a,b,c \in K$, then
$$T \left( \overline{ a \overline{ b \overline{ c \overline{ d\overline{e} } } } } \right)
= T \left( \overline{ \overline{e} \overline{ \overline{d}
\overline{ \overline{c} \overline{ \overline{b} \overline{a} } } } } \right)$$
for all $a,b,c,d,e \in K$.
\end{lem}

\begin{proof}
We iterate (iii) with the aid of (i) and (ii).
Take $a,b,c,d,e \in K$. Then
$$T \left( \overline{ a \overline{ b \overline{
c \overline{ d\overline{e} } } } } \right) =
T \left( \overline{a} \left( b \overline{
c \overline{d \overline{e} } } \right) \right) =
T \left( \left( c \overline{d \overline{e}} \right) \overline{b} \overline{a} \right) =$$
$$=T \left( d\overline{e} \left( \overline{c} \overline{ \overline{b}\overline{a} } \right) \right)
= T \left( e \left( \overline{d} \overline{ \overline{c}
\overline{ \overline{b}\overline{a} } } \right) \right) =
T \left( \overline{ \overline{e} \overline{ \overline{d}
\overline{ \overline{c} \overline{ \overline{b} \overline{a} } } } } \right).$$
\end{proof}

\begin{prop}\label{conwaysmithT}
Suppose that $K$ is left alternative, left distributive,
right weak distributive and has $k$-valued
trace and anisotropic $k$-valued norm.
Suppose also that there is a $k$-linear map
$T : K \rightarrow K$ satisfying the
following three conditions
(i) $T(ab) = T(ba)$;
(ii) $T \left( \overline{a} \right) = T(a)$;
(iii) $T \left( \overline{a(bc)} \right) =
T \left( \overline{c} \left( \overline{b} \overline{a} \right) \right)$
for all $a,b,c \in K$.
If we extend the map $T$ to $K \times^D K$
by the relation $T((a,b)) = T(a)$, for $a,b \in K$,
then (i), (ii) and (iii) hold for all
$a,b,c \in K \times^D K$.
In that case, the norm on $K \times^D K$ is
multiplicative.
\end{prop}

\begin{proof}
(i) and (ii) are straightforward. Now we show (iii). Take
$a,b,c,d,e,f \in K$ and put $x = (a,b)$, $y = (c,d)$ and $z =
(e,f)$. We have to check five different cases.

\vspace{1mm}

\noindent Case 1: at least two of $b$, $d$ and $f$ are equal to
zero. Then, by (iii), we get that
$T \left( \overline{x(yz)} \right) =
T \left( \overline{a(ce)} \right) = T \left( \overline{e}
(\overline{c} \ \overline{a}) \right) = T \left( \overline{z} \left(
\overline{y} \ \overline{x} \right) \right).$

\noindent Case 2: $b = 0$, $d \neq 0$ and $f \neq 0$. Then, by
left distributivity and (iii), we get that
$$T \left( \overline{x(yz)} \right) =
T \left( \overline{(a,0) ( (c,d)(e,f) )} \right) = T \left(
\overline{ a \left( ce + D \overline{ d \overline{f} } \right) }
\right) =$$
$$= T \left( \overline{ a(ce) } \right) +
D T \left( \overline{ a \overline{ d \overline{f} } } \right) = T
\left( \overline{e} (\overline{c} \ \overline{a}) \right) + D T
\left( \overline{f (\overline{d} a)  } \right) =$$
$$= T \left( \overline{e} (\overline{c} \ \overline{a}) +
D \overline{f (\overline{d} a)  } \right) = T \left(
(\overline{e},-f) (\overline{c} \ \overline{a} , \overline{
(\overline{-d})a } ) \right)
%$$= T \left( (\overline{e},-f) ( (\overline{c},-d)(\overline{a},0) ) \right)$$
= T \left( \overline{z} (\overline{y} \ \overline{x}) \right)$$

\noindent Case 3: $b \neq 0$, $d = 0$ and $f \neq 0$. Then,
by left distributivity and (iii), we get that
$$T \left( \overline{x(yz)} \right) =
T \left( \overline{(a,b) ( (c,0)(e,f) )} \right) =
T \left( \overline{ (a,b) (ce,cf) } \right) =$$
$$= T \left( \overline{a(ce)} \right) +
D T \left( b \overline{cf} \right)
= T \left( \overline{e}(\overline{c} \ \overline{a}) \right) + D T
\left( f \overline{\overline{c}b} \right) = T \left(
\overline{e}(\overline{c} \ \overline{a}) + D f
\overline{\overline{c}b} \right) =$$
$$= T \left( (\overline{e},-f) (\overline{c} \ \overline{a} ,
-\overline{c}b) \right) = T \left( \overline{(e,f)} \left(
\overline{(c,0)} \ \overline{(a,b)} \right) \right) = T \left(
\overline{z} (\overline{y} \ \overline{x}) \right)$$

\noindent Case 4: $b \neq 0$, $d \neq 0$ and $f = 0$. Then,
by left distributivity and (iii), we get that
$$T \left( \overline{x(yz)} \right)
= T \left( \overline{ (a,b) ( (c,d) (e,0) ) } \right)
= T \left( \overline{ (a,b) \left( ce , \overline{ \overline{d}
\overline{e} } \right) } \right) =$$
$$= T \left( \overline{ a(ce) } \right) +
D T \left( \overline{b(\overline{d} \overline{e})} \right)
= T \left( \overline{e} (\overline{c} \ \overline{a}) \right) + D T
\left( \overline{e} \overline{d \overline{b}} \right) =$$
$$= T \left( \overline{e} (\overline{c} \ \overline{a}) +
D \overline{e} \overline{d \overline{b}} \right) = T \left(
\overline{e} \left( \overline{c} \ \overline{a} + D \overline{d
\overline{b}} \right) \right)
= T \left( \overline{z} (\overline{y} \ \overline{x}) \right)$$

\noindent Case 5: $b \neq 0$, $d \neq 0$ and $f \neq 0$. Then,
by left distributivity, Lemma \ref{normlemma} and (iii),
we get that
$$  T \left( \overline{x(yz)} \right) = %T \left( \overline{(a,b) ( (c,d)(e,f) )} \right)$$
T \left( \overline{ (a,b) \left( ce + D \overline{ d \overline{f} }
, \overline{\overline{d} \overline{e}} + \overline{ \overline{d}
\overline{ \overline{c} \overline{ \overline{d^{-1}} \ \overline{f}
} } } \right) } \right) =$$
$$= T \left(
\overline{ a(ce + D \overline{ d \overline{f} }) + D b \left(
\overline{d} \overline{e} + \overline{d} \overline{ \overline{c}
\overline{ \overline{d^{-1}} \ \overline{f} } }   \right) }
\right)$$
$$= T \left( \overline{a(ce)} \right) +
D T \left( a \left( \overline{d \overline{f}} \right) \right) + D T
\left( b ( \overline{d} \overline{e} ) \right) + D T \left( b \left(
\overline{d} \overline{ \overline{c} \overline{ \overline{d^{-1}} \
\overline{f} } } \right) \right)$$
$$= T \left( \overline{e} (\overline{c} \ \overline{a} \right) +
D T \left( f ( \overline{d} a ) \right) + D T \left( \overline{e}
\overline{ d \overline{b} } \right) + T \left( f \left( \overline{d}
\overline{ c \overline{ \overline{d^{-1}} \overline{b} } } \right)
\right)$$
$$= T \left( \overline{e} (\overline{c} \ \overline{a} + D \overline{ d \overline{b} }) +
D f \left( \overline{d} \overline{ c \overline{ \overline{d^{-1}}
\overline{b} } } \right) \right)$$
$$= T \left( (\overline{e},-f) \left( \overline{c} \ \overline{a} + D
\overline{ d \overline{b} } , - \overline{ \overline{d}a } -
\overline{ \overline{d} \overline{ c \overline{ \overline{d^{-1}}
\overline{b} } } } \right) \right) = T \left( \overline{z}
(\overline{y} \overline{x}) \right)$$
The last part follows from the above and Propositions
\ref{multiplicativenorm},
\end{proof}

\begin{cor}\label{conwaysmithmultiplicativenorm}
Every $k$-algebra $k^i$, for $i \geq 0$,
as defined in Corollary 3, has multiplicative norm.
\end{cor}

\begin{proof}
Suppose that we define the map $T : k^i \rightarrow k^i$,
for all $i \geq 0$,
inductively by $T = {\rm id}_k$ on $k^0 = k$
and $T((a,b)) = T(a)$ for $(a,b) \in k^{i+1}$.
Since obviously ${\rm id}_k$ satisfies (i), (ii)
and (iii) from Proposition 10, it follows
that the same is true for $T$ on any $k^i$, for $i \geq 0$.
Therefore, by Propositions 4 and 9,
the norm is multiplicative on $k^i$, for $i \geq 0$.
\end{proof}

\end{document}